
%
%
%
%
  %

   \def\pa{\parallel}
   \def\schwm{{\cal S}(M)}
   \def\Z{\bf Z}
   \def\R{\bf R}
   \def\T{\bf T}
   \def\schwr{{\cal S}(\R)}
   \def\GS{{\cal S}(G\times M)}
   \def\GC{{\rm{C}}^{*}(G\times M)}
   \def\GA{{\cal S}(G, A)}
   \def\GB{{\rm{C}}^{*}(G, B)}

  \font\tensc=cmcsc10 
  \font\twelvebf=cmbx10 scaled \magstep1 

  \magnification=\magstep1
  \tolerance=10000
  \hsize=6.5truein
  \vsize=9.2truein


\def\statement#1#2{\medbreak\noindent{\bf#1.\enspace}{\sl#2}\medbreak}
  \def\cite#1{{\rm[\bf #1\rm]}}
  \def\stdbib#1#2#3#4#5#6#7#8{\smallskip \item{[#1]} #2, ``#3'',
    {\sl#4} {\bf#5} (#6), #7--#8.}
  \def\bib#1#2#3#4{\smallskip \item{[#1]} #2, ``#3'', {#4.}}

  \def\BlackadarCuntz {1}
  \def\ConnesThom {2}
  \def\ConnesDiffGeo {3}
  \def\DuCloux {4}
  \def\Palmer {5}
  \def\PhilSchw {6}
  \def\Pytlik {7}
  \def\shortpf {8}
  \def\Dm {9}
  \def\SSI {10}
  \def\rep {11}
  \def\sameauthor{\underbar{\hbox to 1.5truecm{\hfil}}}

  \null
  \vskip -\voffset
  \vskip -1.5truecm
  \baselineskip=9pt

  \vskip 2.5truecm
  \baselineskip=20pt
  \centerline{\twelvebf SUMMARY OF SPECTRAL INVARIANCE RESULTS}

  \bigskip
  \bigskip
  \centerline{\tensc Larry B. Schweitzer}
  \bigskip
  \bigskip

{\narrower\bigskip
\noindent {\bf Abstract.}\   The author's recent results on
spectral invariant dense subalgebras of C*-algebras
associated with dynamical systems are summarized.
If $G$ is a compactly generated polynomial growth
Type R Lie group, and the action of $G$ on $\schwm$ (Schwartz
functions on a locally compact $G$-space $M$) is tempered in
a certain sense, then there is a natural
smooth crossed product ${\cal S}(G \times M)$ which is dense and
spectral invariant in the C*-crossed product ${\rm{C}}^{*}(G \times M)$.
\bigskip\bigskip}
%

The  theory of  differential
geometry on a C*-algebra (or noncommutative space)
Connes \cite{\ConnesDiffGeo} requires
the use of \lq\lq differentiable structures\rq\rq\
for these noncommutative spaces,
or some sort of algebra of \lq\lq differentiable functions\rq\rq\
on the noncommutative space.  Such algebras of functions have usually been
provided by a dense subalgebra of smooth functions $A$
for which the $K$-theory $K_{*}(A)$ is the same as the
$K$-theory of the C*-algebra $K_{*}(B)$
(see for example  Blackadar-Cuntz \cite{\BlackadarCuntz},
and the recent works of  J. Bost, G. Elliott, T. Natsume,
R. Nest,   R. Ji, P. Jollissaint, V. Nistor and many others).

One goal of both of the papers Schweitzer \cite{\SSI} \cite{\rep}
was to construct such dense subalgebras of smooth functions
in the case that $B$ is a C*-crossed product $\GC$
associated with a dynamical system, or more specifically with
an action of a Lie group $G$ (not necessarily connected)
on a locally compact space $M$.
In these papers, we realize this goal by constructing
smooth crossed products
$\GS$ of Schwartz functions on $G\times M$,
which are spectral invariant in  the C*-crossed product.
{\underbar{Spectral}} {\underbar{invariant}}
means that the spectrum of every element
of $\GS$
is the same in $\GS$ and $\GC$.
In the
language of Palmer \cite{\Palmer}, this is the same as saying that
$\GS$ is a {\underbar{spectral}} \underbar{subalgebra} of
$\GC$.
By Schweitzer \cite{\shortpf}, Lemma 1.2, Corollary 2.3, and
Connes \cite{\ConnesThom}, VI.3, spectral
invariant subalgebras have the same $K$-theory as the C*-algebra itself,
so these smooth crossed products $\GS$
provide us with the \lq\lq noncommutative
differentiable structures\rq\rq\  we are looking for.
\vskip\baselineskip

I will begin by describing the results obtained in \cite{\rep}.
The idea in that paper is to employ
the following theorem, which is interesting in its own right.
It gives a condition for a dense Fr\'echet subalgebra $A$ to be a spectral
invariant subalgebra of $B$  when
certain subrepresentations of
{\underbar{topologically}} irreducible representations of $A$ extend
appropriately to $B$.  This is in contrast to the situation in
\cite{\shortpf}, Theorem 1.4, Corollary 1.5,
which says that {\underbar{algebraically}}
{irreducible} representations extend iff the
subalgebra is spectral invariant.
%

If $E$ is an $A$-module, we say that $E$ is \underbar{algebraically}
\underbar{cyclic} iff there exists an $e \in E$ such that the algebraic
span $Ae$ is equal to $E$.
\statement{Theorem 1 (\cite{\rep}, Theorem 1.4)}
{Let $A$ be a dense $m$-convex Fr\'echet
subalgebra
of a C*-algebra $B$ with continuous inclusion map $A\hookrightarrow B$.
Assume that every algebraically cyclic subrepresentation of every
topologically  irreducible representation of $A$ on a Banach space
is contained in a *-representation of $B$ on a Hilbert space.  Then $A$ is
spectral invariant in $B$.}
%

The smooth subalgebras $\GS$ are shown to be
$m$-convex Fr\'echet algebras
in Schweitzer \cite{\Dm}, \S 3, and their topologically irreducible
representations are relatively accessible when the C*-crossed
product is CCR.  Hence Theorem 1 gives many new interesting cases
of spectral invariant smooth crossed products $\GS$.
 For example,
results are obtained
when $G$ is a closed subgroup of a connected, simply connected
nilpotent Lie group with
certain restrictions on the isotropy subgroups (that they be CCR for
one), and when the action of $G$ on $M$ has closed orbits.
(See the examples in \cite{\rep}, \S 18, \S 2, \S 16-17.)

A simple illustrative example is given by $\Z$ acting by translation
on $\R$.  The Schwartz functions $ {\cal S}(\Z \times \R)$
with convolution multiplication provide a
dense subalgebra  of smooth functions of the
C*-crossed product ${\rm{C}}^{*}(\Z \times \R)$.
Any topologically irreducible representation of ${\cal S}(\Z \times \R)$
must factor through an
orbit
to a representation of the convolution algebra
${\cal S}(\Z \times \Z)$ \cite{\rep},
Theorem 14.1.
The latter algebra is
a smooth version of the compact operators, whose
representation theory is quite nice Du Cloux \cite{\DuCloux},
Corollary 3.5
or \cite{\rep}, Theorem 15.1, Example 2.5.  Theorem 1
may then be applied
to obtain the spectral invariance of ${\cal S}(\Z \times \R)$
in ${\rm{C}}^{*}(\Z \times \R)$.

As one might speculate, when the C*-crossed product is not
CCR (or at least when it is not GCR), the representation theory
of the dense subalgebra becomes quite complicated as does the
representation theory of the C*-algebra.  For example, the
dense subalgebra $A_{\theta}^{\infty}$, given by the
canonical
action of $\T^{2}$ on  the irrational
rotation C*-algebra $A_{\theta}$,
is  spectral invariant but does not
satisfy the hypothesis of Theorem 1.  That is, there exists certain
\lq\lq bad\rq\rq \
topologically irreducible representations of $A_{\theta}^{\infty}$,
which have algebraically cyclic subrepresentations which
do not extend to $A_{\theta}$ \cite{\rep}, Example 7.1.
\vskip\baselineskip

In order to
get results in the non-CCR case, a new method is needed to
replace Theorem 1.  Such a method, or methods,
is introduced in \cite{\SSI}, which I shall now describe.

We begin by trying to show that the smooth crossed product $\GA$
is spectral invariant in $L^{1}(G, B)$ instead of in
the C*-crossed product $\GB$.
Let $\pa \quad \pa_{0}$ be the norm on $B$, and let $\bigl\{\pa
\quad \pa_{n}\bigr\}_{n=0}^{\infty}$
be a family of increasing submultiplicative norms giving the
topology of $A$.  In the paper Blackadar-Cuntz \cite{\BlackadarCuntz},
the condition
{$\pa ab \pa_{n} \leq C\sum_{i+j= n} \pa a \pa_{i} \pa b \pa_{j}$}
for all $a, b$ in $A$, is used to show that $A$ is a spectral
invariant subalgebra of $B$.
The commutative Fr\'echet algebra $\schwm$ of Schwartz functions
on $M$ satisfies this condition
in $C_{0}(M)$ \cite{\SSI}, \S 2.
Moreover, for some very nice
actions of $G$ on $A$ (isometric on each norm),
one can  show that
if the norms on $A$ satisfy the condition in $B$,
then the norms on the smooth
crossed product $\GA$ satisfy the condition in
$L^{1}(G, B)$.
The following more general
condition introduced in \cite{\SSI} does the same thing without
requiring an isometric action.

We say that
 $A$ is {\underbar{strongly} \underbar{spectral} \underbar{invariant} }
in $B$ if
$$\eqalign{\quad
 (\exists C>0)(\forall m) & (\exists D_{m}>0)(\exists p_{m}\geq m)
(\forall n)(\forall a_{1},\dots a_{n} \in A)\hfill \cr
& \biggl\{ \pa a_{1} \dots a_{n} \pa_{m}
\leq D_{m}C^{n}
\sum_{k_{1}+ \dots k_{n} \leq p_{m}}
 \pa a_{1} \pa_{k_{1}}\dots  \pa a_{n} \pa_{k_{n}}
\biggr\}.\cr}\leqno (*)  $$
Notice that in the summand of (*), at most $p_{m}$ of the natural numbers
$k_{j}$
are nonzero, regardless of $n$.
The idea behind showing that strong spectral invariance implies
spectral invariance is given by setting $a_{1}= \dots a_{n} = a$ in (*).
We have
$$ \eqalign{ \pa a^{n} \pa_{m} & \leq  DC^{n} \sum_{k_{1} + \dots k_{n}
\leq p}  \pa a \pa_{k_{1}} \dots \pa a \pa_{k_{n}} \cr
& \leq K^{n}  \pa a \pa_{0}^{n-p}  \pa a \pa_{p}^{p}\cr}
$$
where  $p$ is fixed as $n$ runs.
It follows that
the series $(1-a)^{-1} = 1 + a + a^{2} + \dots $ converges absolutely in
the norm $\pa \quad \pa_{m}$ when $\pa a \pa_{0}$ is sufficiently small.
So $1-a$ is invertible in the completion of $A$ in $\pa \quad \pa_{m}$
when $a$ is sufficiently close to $0$ {\underbar{in} \underbar{$B$}}.
The rest of the argument
is in Theorem 1.17 of \cite{\SSI}.

There are also examples of spectral invariant dense subalgebras
which are not strongly spectral invariant \cite{\SSI}, Example 1.13.
The following theorem and corollary illustrates the usefulness of
the concept of strong spectral invariance.

We say that a Lie group $G$ (not necessarily connected) is
{\underbar{compactly} \underbar{generated}}
if $G$ has an open relatively compact
neighborhood $U$ of the identity which satisfies
$\bigcup_{n=0}^{\infty}
U^{n} = G$ and $U^{-1}= U$.
We call $\tau(g) =\min \bigl\{ \, n\,\bigl| \, g \in U^{n} \,\bigr\}$
the {\underbar{word}} \underbar{gauge} on $G$.
(The smooth crossed product $\GA$ is then defined to be
the set of $G$-differentiable $\tau$-rapidly vanishing
functions from $G$ to $A$.)
We say that the action of $G$ on $A$ is \underbar{$\tau$-
tempered} if for
every $m$, $\pa \alpha_{g}(a) \pa_{m} $ is bounded by a
polynomial in $\tau(g)$ times $\pa a \pa_{n}$ for some $n$.
Finally, we say that  $G$
is  {\underbar{Type R}} if
 all the eigenvalues of $Ad_{g}$ lie on the unit
circle for each $g \in G$.
\statement{Theorem 2}{If $A$ is strongly spectral invariant in $B$
and $G$ is a compactly generated Type R  Lie group,
and the action of $G$ on $A$ is $\tau$-tempered,
then the
smooth crossed product $\GA$ is strongly spectral invariant
in $L^{1}(G, B)$.}
\statement{Corollary 3}{For compactly generated Type R  Lie groups,
for which the action of $G$ on $\schwm$ is $\tau$-tempered,
the smooth crossed product $\GS$
is strongly spectral invariant in $L^{1}(G, C_{0}(M))$.}
%
\noindent
(Note that $\GS$ is shorthand for ${\cal S}(G, \schwm)$.)
It is the subadditivity of $\tau$ and the strong spectral invariance
of $A$ in $B$
that play the essential role in the proof of Theorem 2 and Corollary 3.
The hypotheses that $G$ be Type R and that the action is $\tau$-tempered
are not used in the proof,
but they are necessary to assure the
existence of the smooth crossed product $\GS$,
and to assure that $\GS$ is a Fr\'echet *-algebra.
There are a wide variety of Type R Lie groups
(see below or \cite{\Dm}, \S 1.4,
\cite{\SSI}), and also many examples of $\tau$-tempered actions of
such groups $G$ on $\schwm$ \cite{\Dm}, \S 5, \cite{\SSI}, Examples  6.26-7,
7.20, \cite{\rep}.

{\bf Remark.\ } We clarify what the $\tau$-tempered assumption
can mean in practice.  Let $G$ be the integers $\Z$, and let $G$ act
on $\R$ via $\alpha_{n}(r) = e^{-n}r$.  The word gauge $\tau$ is
equivalent in an appropriate sense to the absolute value function
$\tau(n) = |n|$.  If we take $\schwm= C_{0}(\R)$, or
$\schwm = C_{0}^{\infty}(\R)$, then $\alpha$ is an isometric
action of $\Z$ on $\schwm$, meaning that $\alpha$
leaves each seminorm invariant, and so $\alpha$ is  $\tau$-tempered.
However, if we take $\schwm$ to be the
standard Schwartz functions $\schwr$,
then for a fixed  $\varphi\in \schwm$,
$\pa \alpha_{n}(\varphi) \pa_{m}$
will in general grow
exponentially in $n$ as $n \rightarrow +\infty$, so $\alpha$
is no longer $\tau$-tempered.
So the
$\tau$-temperedness condition does place a restriction on what
$\schwm$ can be for a given action.
If the action of $\Z$ on $\R$ were by translation $\alpha_{n}(r)=
r + n$, then the action of $\Z$ on $\schwr$ {\underbar{would}}
be $\tau$-tempered.
For general $M$ and $G$ as in Corollary 3, and
regardless of the action,
one can always get a $\tau$-tempered action
by taking $\schwm= C_{0}(M)$, or $\schwm= C_{0}^{\infty}(M)$,
where the superscript $\infty$ means \lq\lq $G$-differentiable\rq\rq.
\vskip\baselineskip
%

Note that Corollary 3 makes
no assumption about the crossed product being CCR. No restrictions
on the action of $G$ or
the isotropy subgroups are needed.
However,
we are left with the question of whether $\GS$ is
spectral invariant in the C*-crossed product $\GC$,
and not just $L^{1}(G, C_{0}(M))$.
To take care of this we generalize a result of Pytlik \cite{\Pytlik}
which says that if $G$ has  \underbar{polynomial} \underbar{growth}, then
the rapidly vanishing $L^{1}$-functions on $G$ form a  symmetric
Fr\'echet *-algebra, which consequently is spectral invariant in
$C^* (G)$.  In particular, we show in \S 7 of \cite{\SSI} that the
rapidly vanishing $L^{1}$-functions from $G$ to $B$
is  spectral invariant in the C*-crossed product $\GB$ when
$G$ has polynomial growth.   Since these rapidly vanishing
$L^{1}$-functions are also spectral
invariant in $L^{1}(G, B)$, and since they contain
the smooth crossed product,
we are able to conclude that
the smooth crossed product is spectral invariant in the
C*-crossed product when $G$ has polynomial growth.
Our main result is then:
\statement{Corollary 4}{For compactly generated polynomial growth
Type R Lie groups $G$, and $\tau$-tempered actions of $G$ on
$\schwm$,
the smooth crossed product $\GS$
is spectral invariant in the C*-crossed product $\GC$.}
Examples of such groups are given by
finitely generated polynomial growth discrete groups,
 compact or  connected nilpotent
Lie groups,  the group of Euclidean motions on the plane,
any motion group, or any closed  subgroup of one of these.
Numerous examples of smooth crossed products which are spectral invariant
because of  Corollary 4 can be found in \cite{\SSI}, Examples 2.6-7,
6.26-7, 7.20, \cite{\rep}, \cite{\Dm}, \S 5.

We remark that in \cite{\PhilSchw}, methods are given to show that
$\GS\hookrightarrow \GC$ is an isomorphism on $K$-theory
without using spectral invariance,
whenever $G$ is a closed subgroup of a connected, simply connected
nilpotent Lie group, and the action of $G$ on $\schwm$ is $\tau$-tempered
\cite{\PhilSchw}, Example 3.2.

  \bigbreak
  \centerline{\tensc References}
  \nobreak\medskip
  \frenchspacing

  \bib{\BlackadarCuntz}
  {B. Blackadar}
  {Differential Banach algebra norms and smooth subalgebras of C*-algebras}
  {\it J. Operator Theory} {to appear.}
  \stdbib{\ConnesThom}
  {A. Connes}
  {An analogue of the Thom
isomorphism for crossed products of a C*-algebra by an action of $\R$}
  {Adv. in Math.} {39} {1981} {31} {55}
  \stdbib{\ConnesDiffGeo}
  {\sameauthor}
  {Non-commutative differential geometry}
  {Publ. Math. I.H.E.S.} {62} {1985} {257} {360}
  \stdbib{\DuCloux}
  {F. Du Cloux}
  {Repr\'esentations temp\'er\'ees des groupes de Lie nilpotent}
  {J. Funct. Anal.} {85} {1989} {420} {457}
  \stdbib{\Palmer}
  {T. W. Palmer}
  {Spectral Algebras}
  {Rocky Mountain J. Math.} {22(1)} {1992} {293} {328}
  \bib{\PhilSchw}
  {N. C. Phillips and L. B. Schweitzer}
  {Representable $K$-theory of smooth crossed products by $\R$ and $\Z$}
  {preprint} {(1992).}
  \stdbib{\Pytlik}
  {T. Pytlik}
  {On the spectral radius of elements in group algebras}
  {Bulletin de l'Acad\'emie Polonaise des Sciences, astronomiques et
physiques}  {21} {1973} {899} {902}
  \stdbib{\shortpf}
  {L. B. Schweitzer}
  {A short proof that $M_{n}(A)$ is local if $A$ is local and Fr\'echet}
  {Intl. Jour. Math.} {3(4)} {1992} {581} {589}
  \bib{\Dm}
  {\sameauthor}
  {Dense $m$-convex Fr\'echet subalgebras of operator algebra
crossed products by Lie groups}
  {preprint} {(1992).}
  \bib{\SSI}
  {\sameauthor}
  {Spectral invariance of dense subalgebras of operator algebras}
  {\it Intl. Jour. Math} {to appear.}
  \bib{\rep}
  {\sameauthor}
  {Representations of dense subalgebras of C*-algebras}
  {preliminary version} {(January, 1993).}

 \bigskip
 \baselineskip=12pt
 \line{\it Department of Mathematics and Statistics, University of
Victoria,\hfill}
 \line{\it Victoria, B.C., Canada V8W 3P4. \ \ E-mail:
lschweit@sol.uvic.ca\hfill}

\end